\documentclass[11pt]{article}
\usepackage{amsmath, amssymb, amsthm}
\usepackage[all]{xypic}

\theoremstyle{definition}

\theoremstyle{remark}

\newtheorem{ntn}[equation]{Notation}

\makeatletter
\renewcommand{\subsection}{\@startsection{subsection}{2}{0pt}{-3ex
plus -1ex minus -0.2ex}{-2mm plus -0pt minus
-2pt}{\normalfont\bfseries}} \makeatother

\numberwithin{equation}{subsection}

\newdir^{ (}{{}*!/-5pt/@^{(}}

 \newcommand{\iso}{{\;\stackrel{_\sim}{\to}\;}}
\newcommand{\beq}{\begin{equation}\label}
\newcommand{\eeq}{\end{equation}}
\def\ccirc{{{}_{\,{}^{^\circ}}}}

\newcommand{\tooo}{{\;{-\!\!\!-\!\!\!-\!\!\!-\!\!\!\longrightarrow}\;}}

\newcommand{\pr}{\operatorname{pr}}

\def\k{\mathbf{k}}

\def\lb{{\mathbf{l}}}
\def\ppi{\mathrm{pr}}
\def\o{\otimes}
\def\h{\mathrm{h}}

\def\dq{\overline{Q}}

\def\Id{\mathrm{Id}}

\def\moy{{\mathrm{Moyal}}}
\def\Rep{{\mathrm{Rep}}}
\def\d{{\mathbf{l}}}
\def\Sym{{\text{Sym\ }}}
\def\Z{{\mathbb Z}}
\def\tr{{\text{tr}}}

\begin{document}

\centerline{\Large{\textbf{Moyal quantization of necklace Lie algebras}}}

\bigskip

\centerline{\sc{Victor Ginzburg and Travis Schedler}}
\medskip

\begin{abstract} We use Moyal-type formulas to
construct a Hopf algebra quantization of the
necklace Lie bialgebra associated with a quiver.
\end{abstract}

\section{Introduction}
\subsection{Reminder on Moyal product.}\label{reminder} Let $V$ be a finite
dimensional vector space equipped with a nondegenerate
bivector $\pi\in\wedge^2V$. Associated with $\pi$ is
a Poisson bracket $f,g\mapsto \{f,g\}:=\langle df\wedge dg,\pi\rangle$
on $\k[V],$ the polynomial algebra on $V$.
The usual commutative product $m: \k[V]\o \k[V]\to\k[V]$
and the Poisson bracket $\{-,-\}$ make $\k[V]$ a  Poisson algebra.
This Poisson algebra has a well-known {\em Moyal-Weyl quantization} 
(\cite{M}, see also \cite{CP}).
This is
an associative   star-product 
depending on a  formal quantization parameter $\h$, defined by the formula
\beq{star}
f *_\h g:=m\ccirc e^{\frac{1}{2} \h \pi} (f\o g)\in \k[V][\h],\quad
\forall f,g\in\k[V][\h].
\eeq

To explain the meaning of this formula, view elements of $\Sym V$
as
constant-coefficient differential operators on $V$. Hence,
an element of $\Sym V\o \Sym V$ acts  as a
constant-coefficient differential operator on 
the algebra $\k[V]\o\k[V]=\k[V\times V].$ Now, identify
$\wedge^2V$ with the subspace of skew-symmetric tensors in $V\o V$.
This way, the bivector $\pi\in\wedge^2V\subset V\o V$ becomes
a second order constant-coefficient differential operator
$\pi: \k[V]\o\k[V]\to\k[V]\o\k[V].$ Further, it is clear that
for any element $f\o g\in\k[V]\o\k[V]$ of total degree
$\leq N$, all terms with $d>N$ in the
infinite sum $e^{\h\cdot \pi}(f\o g)=\sum_{d=0}^\infty
\frac{\h^d}{d!} \pi^d(f\o g)$
vanish, so the sum makes sense.
Thus,
the 
symbol $m\ccirc e^{\h\cdot \pi}$ in the right-hand side
of formula \eqref{star} stands for the composition
$$ \k[V]\o\k[V]\stackrel{e^{\h\cdot \pi}}\tooo \k[V]\o\k[V]\o\k[\h]
\stackrel{m\o\Id_{\k[\h]}}\tooo  \k[V]\o\k[\h],
$$
where
$e^{\h\cdot \pi}$
is an infinite-order formal differential operator.

In down-to-earth terms, choose coordinates
$x_1, \ldots, x_n, y_1, \ldots, y_n$ on $V$ such that
the bivector $\pi$, resp., the Poisson bracket $\{-,-\}$, takes 
the canonical form
\beq{pois}
\pi = \sum_i \frac{\partial}{\partial x_i} \o \frac{\partial}{\partial
y_i} - \frac{\partial}{\partial y_i} \o \frac{\partial}{\partial x_i},
\quad\text{resp.,}\quad
\{f,g\}=\sum_i \frac{\partial f}{\partial x_i}\frac{\partial g}{\partial
y_i} - \frac{\partial  f}{\partial y_i} \frac{\partial g}{\partial x_i}.
\eeq
Thus, in canonical coordinates
$x=(x_1, \ldots, x_n), y=(y_1, \ldots, y_n),$ formula \eqref{star} for the Moyal
product
reads
\begin{align}\label{star1}
(f *_\h g)(x,y)&=\sum_{d=0}^\infty
\frac{\h^d}{d!}\left(
\sum_i \frac{\partial}{\partial
x'_i} \frac{\partial}{\partial
y''_i} - \frac{\partial}{\partial y'_i} 
\frac{\partial}{\partial x''_i}
\right)^df(x',y') g(x'',y'')\Big|_{{x'=x''=x}\atop{y'=y''=y}}\nonumber\\
&=\sum_{\mathbf{j},\mathbf{l}\in\Z^n_{\geq 0}}
(-1)^{\mathbf{l}|}\frac{\h^{|\mathbf{j}|+|\mathbf{l}|}}{\mathbf{j}!\,\mathbf{l}!}\cdot
\frac{\partial^{\mathbf{j}+\mathbf{l}}f(x,y)}{\partial x^\mathbf{j}\partial y^\mathbf{l}}
\cdot
\frac{\partial^{\mathbf{j}+\mathbf{l}}g(x,y)}{\partial y^\mathbf{j}\partial x^\mathbf{l}},
\end{align}
where for $\mathbf{j}=(j_1, \ldots,j_n)\in \Z^n_{\geq 0}$
we put $|\mathbf{j}|=\sum_i j_i$ and given
$\mathbf{j},\mathbf{l}\in \Z^n_{\geq 0},$ write
$$\frac{1}{\mathbf{j}!\,\mathbf{l}!}\frac{\partial^{\mathbf{j}+\mathbf{l}}}
{\partial x^{\mathbf{j}}\partial y^{\mathbf{l}}}:=
\frac{1}{j_1!\ldots
j_n!l_1!\ldots l_n!}\cdot\frac{\partial^{|\mathbf{j}|+|\mathbf{l}|}}
{\partial x_1^{j_1}\ldots\partial x_n^{j_n}\partial
y_1^{l_1}\ldots\partial y_n^{l_n}}.
$$

A more conceptual approach to formulas \eqref{star}--\eqref{star1}
is obtained by introducing the {\em Weyl algebra} $A_\h(V)$.
This is a $\k[\h]$-algebra defined by the quotient
$$
A_\h(V):= (TV^*)[\h]/I(u\o u' - u'\o u-\h\cdot\langle \pi,
u\o u'\rangle)_{u,u'\in V^*},
$$
where $TV^*$ denotes the tensor algebra of the vector space $V^*$,
and $I(\ldots)$ denotes the two-sided ideal generated by the
indicated set.
Now, a version of the Poincar\'e-Birkhoff-Witt theorem
says that the natural {\em symmetrization map}
yields a  $\k[\h]$-linear bijection
$\phi_W: \k[V][\h]\iso A_\h(V)$.
Thus, transporting the multiplication map in the Weyl algebra $A_\h(V)$ via
this bijection, one obtains an associative product
$$\k[V][\h]\o_{\k[\h]}\k[V][\h]\to\k[V][\h],
\quad f\o g\mapsto \phi_W^{-1}(\phi_W(f)\cdot \phi_W(g)).
$$
It is known that this  associative product is equal to the one
given by formulas  \eqref{star}--\eqref{star1}.

\subsection{The quiver analogue.} The goal of this paper 
is to extend the constructions outlined above to noncommutative
symplectic geometry. Specifically, following
an original idea of Kontsevich \cite{K}, to any quiver, one
 associates canonically a certain Poisson algebra 
(\cite{G},  \cite{BLB}).  Then, 
we will produce a quantization of that Poisson algebra
given by an explicit formula analogous to 
 formulas  \eqref{star}--\eqref{star1}.

In more detail, fix a quiver with vertex set $I$ and edge set $Q,$ and
let $\dq$ be the double of $Q$ obtained by adding
reverse edge $e^*\in\dq$ for each edge $e\in Q$.
Let $P$ be the {\em path algebra} of $\dq$. The commutator
quotient space $P/[P,P]$ may be identified naturally
with the space $L$ spanned by cyclic paths (forgetting which was
the initial edge), sometimes called
{\em necklaces}. Letting $\text{pr}_L: P \rightarrow P/[P,P] = L$ be he
projection, there is a natural bilinear
pairing 
\beq{pair}\{-,-\}:\
L\o L\to L,\quad
f\o g\mapsto\{f,g\}:= \text{pr}_L \biggl( \sum_{e \in Q} \frac{\partial f}{\partial e} 
\frac{\partial g}{\partial e^*} - \frac{\partial f}{\partial e^*} 
\frac{\partial g}{\partial e} \biggr).
\eeq
Interpreting $\frac{\partial}{\partial e}, \frac{\partial}{\partial e^*}$
appropriately as maps $L \rightarrow P, P \rightarrow P$, this formula,
which is a quiver analogue of \eqref{pois}, provides $L$
with a Lie algebra structure, first studied in \cite{G},  \cite{BLB}.
More recently, the second author showed in \cite{S}
that there is also a natural Lie {\em cobracket}
on $L$.
To explain this, write $a_1\cdots a_p\in P$ for a path
of length $p$ and let
 $1_i$ denote the trivial (idempotent) path at the
vertex $i\in I$. Further, for
 any edge $e\in \dq$ with head $h(e)\in I$ and tail 
$t(e)\in I$, 
 let
$D_e: P\to P\o P$ be
the derivation defined by the assignment
$$D_e:\
P\to P\o P,\quad
a_1\cdots a_p\mapsto\sum_{a_r=e} 
a_1\cdots a_{r-1}1_{t(e)}\o1_{h(e)} a_{r+1}\cdots a_p.
$$
The map $D_e$ is a derivation. Moreover, the
following map,  cf. \cite[(1.7)-(1.8)]{S}:
\beq{delta}
\delta: L\to L\wedge L,
\quad
f\mapsto \delta(f)= (\text{pr}_L \o \text{pr}_L) 
\biggl( \sum_{e \in Q} D_e(\frac{\partial f}{\partial e^*})
- D_{e^*}(\frac{\partial f}{\partial e}) \biggr)
\eeq
(that is, in a sense, dual to \eqref{pair})
makes the Lie algebra $L$ a Lie {\em bialgebra},
to be referred to as the {\em necklace Lie bialgebra}.

The necklace Lie bialgebra admits a very interesting quantization.
Specifically, the main construction of \cite{S} produces
 a Hopf $\k[\h]$-algebra
$A_\h(Q)$ equipped with an algebra
isomorphism $A_\h(Q)/\h\cdot A_\h(Q)\iso \Sym L,\,f\mapsto\pr{f}.$
The algebra $A_\h(Q)$ is a quantization of the
Lie bialgebra $L$ in the sense that $A_\h(Q)$ is flat over $\k[\h]$ and, 
for any $a,b\in A_\h(Q),$ one has
$$
\pr\left(\frac{ab-ba}{\h}\right)=\{\pr{a},\pr{b}\},\quad
\text{and}\quad
\pr\left(\frac{\Delta(a)-\Delta^{op}(a)}{\h}\right)=\delta(\pr(a)),
$$
where $\Delta: A_\h(Q)\to A_\h(Q)\o_{\k[\h]}A_\h(Q)$ denotes the
coproduct in the Hopf algebra $A_\h(Q)$,
and where $\Delta^{op}$ stands for the map
$\Delta$ composed with the flip of the two factors in
$A_\h(Q)\o_{\k[\h]}A_\h(Q).$

\subsection{Moyal quantization for quivers.} In  \cite{S}, the Hopf algebra
$A_\h(Q)$ was defined, essentially, by generators
and relations. Thus, the algebra
$A_\h(Q)$  may be viewed, roughly, as a quiver analog 
of the Weyl algebra $A_\h(V)$. One of the main
results proved in  \cite{S} is a version of
Poincar\'e-Birkhoff-Witt (PBW) theorem. The PBW theorem
 insures 
that $A_\h(Q)$ is isomorphic
to $(\Sym L)[\h]$ as a $\k[\h]$-module, in particular, it is
flat over $\k[\h]$. 

The goal of the present paper is to
provide an alternative construction of the Hopf algebra
$A_\h(Q)$. Instead of defining the algebra
by generators and relations, we define
a  multiplication $m$ and comultiplication
$\Delta$ on the vector space
$(\Sym L)[\h]$ by explicit formulas which
are both analogous to 
formula \eqref{star} for the Moyal star-product. In fact, suitably interpreted,
they will be written as $f *_\h g = m \circ e^{\frac{1}{2} \h \pi}(f \o g)$
and $\Delta_h(f) = e^{\frac{1}{2} \h \pi} f$.
We directly check 
 associativity, coassociativity
and compatibility of $m$ and $\Delta$.
Thus, the present approach is (up to some difficulties
involving the antipode) independent of that used in \cite{S}.

Further, in complete analogy with the case of Moyal-Weyl
quantization, we construct a symmetrization
map $\Phi: (\Sym L)[\h]\to A_\h(Q)$. This map is
a bijection, and we show that Hopf algebra structure
on $(\Sym L)[\h]$ defined in this paper may
be obtained by transporting the Hopf algebra structure
on $A_\h(Q)$ defined in \cite{S} via $\Phi$.  

\subsection{Representations for the Moyal quantization.}
In \cite{G}, an interesting representation of the necklace Lie algebra is presented
which is quantized in \cite{S}.  Namely, for any representation of the
double quiver $\overline{Q}$ assigning to each arrow $e \in \overline{Q}$ the
matrix $M_e: V_{t(e)} \rightarrow V_{h(e)}$, we can consider the map $L \rightarrow
\k$ given by $e_1 e_2 \cdots e_m \mapsto \tr(M_{e_1} M_{e_2} \cdots M_{e_m})$.  More
generally, if $\mathbf{l} \in \Z_{\geq 0}^I$, then we can consider the representation
space $\mathrm{Rep}_{\mathbf{l}}(\dq)$ of representations with dimension vector $\mathbf l$,
meaning that $\mathrm{dim}\ V_i = l_i$.  Then this is a vector space of dimension
$\sum_{e \in \dq} l_{t(e)} l_{h(e)}$.  It has a natural bivector $\pi((M_{e})_{ij},
(M_{f})_{kl}) = \delta_{il} \delta_{jk} [e,f]$, where $[e,f] = 1$ if $e \in Q, f = e^*$
and $[e,f] = -1$ if $f \in Q, e = f^*$, with $[e,f]=0$ otherwise. We then have 
the Poisson algebra homomorphism
\beq{trrep}
\tr_{\mathbf{l}}: \Sym L \rightarrow \k[\mathrm{Rep}_{\mathbf{l}}(\dq)], \quad 
\tr_{\mathbf{l}}(e_1 e_2 \cdots e_m)(\psi) = \tr(M_{e_1} M_{e_2} \cdots M_{e_m}).
\eeq

In \cite{S}, this representation was quantized by a representation $\rho_\lb: A \rightarrow
\mathcal D(\mathrm{Rep}_\lb(Q))$, where the latter is the space of differential operators
with polynomial coefficients on $\mathrm{Rep}_\lb(Q)$.  We may modify the $\rho_\lb$ 
and $A$ slightly to obtain $\rho_\lb^\h, A_\h$ so that we have the following diagram:
\begin{equation}
\xymatrix{
\Sym L \ar[rr]_-{\mathrm{asympt. inj.}}^{\tr_\d} & & \k[\Rep_\d(\dq)] \\
A_\h \ar[u] \ar[rr]_-{\mathrm{asympt. inj.}}^{\rho_\d^\h} & &
D_Q
\ar[u]}
\end{equation}
Here, $A_\h$ is obtained from $A$ by modifying (3.3) in \cite{S} so that the right-hand
side has an $\h$ just like (3.4). [Note: More generally, it makes sense to consider
the space where (3.3) has an independent formal parameter $\hbar$; for the Moyal
version, we want the two to be the same.] Then, the representations $\rho_\lb^\h$
send elements $(e_1, 1) (e_2, 2) \cdots (e_m, m) \in A_\h$ (see \cite{S}: this is one
lift of $e_1 e_2 \cdots e_m \in L$)
to operators $\sum_{i_1, i_2, \cdots, i_m} 
\iota(e_1)_{i_1 i_2} \iota (e_2)_{i_2 i_3} \cdots \iota(e_m)_{i_m i_1}$, 
where $\iota(e)$ is the matrix $M_e$ if $e \in Q$, and $\iota(e^*) = M_{e^*}$ for
$e \in Q$, where $M_{e^*}$ is the matrix given by 
$(M_{e^*})_{ij} = -\h \frac{\partial}{\partial (M_e)_{ji}}$.  Then, the space $D_Q \subset
\mathcal D(\mathrm{Rep}_\lb(Q))$ is just generated by $e_{ij}, 
-\h \frac{\partial}{\partial e_{ji}}$.

The diagram indicates that the representations are ``asymptotically injective'' in
the sense that the kernels of the representations $\rho_\lb, \tr_\lb$ have zero
intersection, and moreover, for any finite-dimensional vector subspace $W$ of the
algebra $A$, there is a vector $\lb \in N^I$ such that for each
$\lb' \geq \lb$ (i.e.~such that $l_i' \geq l_i, \forall i$, 
we have that $W \cap \text{Ker }\tr_\lb = 0$ (and similarly for $\rho$).

By construction of the map $\Phi_W$,
the Moyal quantization fits into a diagram as follows:
\begin{equation} \label{2d}
\xymatrix{ \Sym L[\h]_\moy \ar[rr]_{\mathrm{asympt. inj.}}^{\tr_\lb[\h]} \ar[d]
\ar@/_5pc/[dd]^{\Phi_W}_{\sim} & & \k[\h][\Rep_\lb(\dq)]_\moy \ar[d]
\ar@/^5pc/[dd]^{\phi_W}_{\sim} \\ \Sym L
\ar[rr]_-{\mathrm{asympt. inj.}}^{\tr_\lb} & & \k[\Rep_\lb(\dq)] \\ A_\h \ar[u]
\ar[rr]_-{\mathrm{asympt. inj.}}^{\rho^h_\lb} & & D_Q \ar[u] }
\end{equation}
Here, we denote by $\k[\h][\Rep_\lb(\dq)]_\moy$ the Moyal quantization
of $\k[\Rep_\lb(\dq)]$ using the bivector $\pi$, and by $\Sym
L[\h]_\moy$ the quiver version to be defined in this article.  Because
of the asymptotic injectivity, to prove that a Moyal quantization
exists completing the diagram, all that is necessary is the map
$\Phi_W$; then the definitions of the product, coproduct, and antipode
follow.  However, the definitions are interesting in their own right.

\subsection{Organization of the article.}
The article is organized as follows: In Section \ref{mps}, we
will define the Moyal product $*_\h$ on $\Sym L[\h]$. In Section \ref{phiws},
we define the map $\Phi_W$.  Next, in Section \ref{mpts}, we show 
that this transports the product on
$A_\h$ to the product $*_\h$. Finally, in Section \ref{ass}, we directly prove the
associativity of $*_\h$.

In Section \ref{cps} we define
the Moyal coproduct $\Delta_\h$.  Then, in Section \ref{cpts}, we show that $\Delta_\h$
is obtained by transporting the coproduct from $A_\h$ using $\Phi_W$. 
Section \ref{casss} proves directly that $\Delta_\h$ is coassociative, and
Section \ref{bas} shows directly that $*_\h, \Delta_\h$ are compatible, inducing
a bialgebra structure on $\Sym L[h]_\moy$.

In Section \ref{as} we give the definition
of antipode $S$, which clearly is the one obtained from $A_\h$ by transportation.
This makes $\Sym L[h]_\moy$ a Hopf algebra satisfying $S^2 = \Id$.  The eigenvectors
of $S$ are just products of necklaces, with eigenvalue $\pm 1$ depending on the parity
of the number of necklaces.

We will make use of the following tensor convention throughout:
\begin{ntn}
If $S, T$ are $\k[\h]$-modules, then we will always mean by $S \o T$ the
tensor product over $\k[\h]$ (never over $\k$).
\end{ntn}

\subsection{Acknowledgements}
The authors would like to thank Pavel Etingof for useful discussions.
The work of both authors was partially supported by the NSF.

\section{The Moyal product}
\subsection{Definition of the Moyal product $*_\h$.}\label{mps}
To define the product $*_\h$ on $\Sym L[\h]_\moy$, we proceed
by analogy: let $\pi = \sum_{e \in Q} \frac{\partial}{\partial e} \o
\frac{\partial}{\partial e^*} - \frac{\partial}{\partial e^*} \o
\frac{\partial}{\partial e}$.  For each $n \geq 0$, we define an
operator $\pi^n: \Sym L \o \Sym L \rightarrow \Sym L$, and hence
$e^{\frac{1}{2} \h \pi}: \Sym L[\h] \o \Sym L[\h] \rightarrow \Sym L[\h]$ as
follows.  We define the action of each
\begin{equation} \label{act}
T = \frac{\partial}{\partial a_1} \frac{\partial}{\partial a_2}
\cdots \frac{\partial}{\partial a_m} \o \frac{\partial}{\partial
a_1^*} \frac{\partial}{\partial a_2^*} \cdots \frac{\partial}{\partial
a_m^*}, \quad a_i \in \dq, (e^*)^* := e;
\end{equation}
and extend by linearity.  This action is best described by considering
monomials in $\Sym L$ to be collections of closed paths in $\dq$.
Each closed path corresponds to a single cyclic monomial of $L$, so a
collection of closed paths corresponds to a symmetric product of the
corresponding cyclic monomials, giving an element of $\Sym L$.  Such
elements generate all of $\Sym L$.

Take any operator of the form \eqref{act}, and two elements $P, R \in
\Sym L$, which are symmetric products (i.e.~collections) of closed
paths.  Then the element $T$ of \eqref{act} acts on $P \o R$ by
summing over all ordered choices of distinct instances of edges $e_1,
e_2, \cdots, e_m$ in the graph of $P$ such that $e_i$ is identical
with $a_i$ as elements of $\dq$, and over all ordered choices of
distinct instances of edges $f_1, f_2, \cdots, f_m$ in the graph of
$R$ such that $f_i$ is identical with $e_i^*$ as elements of $\dq$,
and adding the following element: Delete each $e_i$ from $P$ and each
$f_i$ from $R$, and join $P$ and $R$ at each $h(e_i) = t(f_i)$ and
each $h(f_i) = t(e_i)$.  The result is some element $Z \in \Sym L$
obtained from $P \o R$, which is some new collection of closed paths (or
isolated vertices, which correspond to idempotents).
So, $T(P \o R)$ is the sum of all such elements $Z$ (some of them can
be identical, of course; we are summing over the element $Z$ we get
for each choice of instances of the given edges in $P$ and $R$).

Let us more precisely define this deletion and gluing process (as
in \cite{S}).  We can
define an ``abstract edge'' of an element
\begin{multline} \label{pform}
P = a_{11} a_{12} \cdots a_{1 l_1} \& a_{21} a_{22} \cdots a_{2 l_2}
\& \cdots \& a_{k1} a_{k2} \cdots a_{k l_k} \\  \& v_1 \& v_2 \& \cdots
\& v_q \in \Sym L
\end{multline}
to be an index $(i,j)$ where $1 \leq i \leq k$ and $1 \leq j \leq
l_k$.  Here we note that $v_i \in I$, the set of vertices of the
quiver, which act as idempotents in the path algebra of the double
quiver.  These indices are considered as edges, just where we keep
track of which occurrence of the edge of $\dq$ we are considering.  Let
$X$ be the set of abstract edges of such an element $P$; then there is
a natural map $\ppi_X: X \rightarrow \dq$ which gives the element of
$\dq$ defined by the given edge.

To cut and glue for a single such element $P$, we need a set of
``cutting edges,'' $I \subset X$, along with a (fixed-point free)
involution $\phi: I \rightarrow I$ such that $\ppi_X \circ \phi = *
\circ \ppi_X$, where $*: \dq \rightarrow \dq$ is the edge reversal
operation.  Then we can define a map $f: X \rightarrow X$ which takes
each edge $(i,j) \notin I$ to the next edge, $(i,j+1)$ (where $j+1$ is
taken modulo $l_i$); and takes each edge $(i,j) \in I$ to $\phi(i,j) +
1$, where the ``$+1$'' operation is just $(i,j) + 1 = (i,j+1)$, again
where $j+1$ is taken mod $l_i$.  The map $f$ is bijective, and
each orbit of $X$ under $f$ is of the form $(x_1, x_2, \ldots, x_p)$
where $f(x_i) = x_{i+1}$, taken modulo $p$.  Each such orbit defines a
cyclic monomial or idempotent as follows: for each $x_i$, let
$\ppi'(x_i) = \ppi_X(x_i)$ if $x_i \notin I$, and $\ppi'(x_i) = t(x_i)$,
the starting vertex idempotent, if $x_i \in I$.  So $\ppi'$ extends
to $\ppi'(x_1, x_2, \ldots, x_p) = \ppi'(x_1) \ppi'(x_2) \cdots \ppi'(x_p)
\in L$, which gives us the desired cyclic monomial or vertex
idempotent.  Then the result of cutting and gluing along the edges $I$
is simply the symmetric product of $\ppi'$ applied to all orbits of $X$
under $f$, symmetric-multiplied by $v_1 \& v_2 \& \cdots \& v_q$ (the
original vertex idempotents are ``untouched'' by cutting and gluing at
edges).

Given two elements $P, R$ of the form \eqref{pform} (except for
different indices $i_l, k, m$, and different edges $a_{ij}$ etc.), we
can cut and glue $P$ and $R$ in an analogous way as follows: Let $X,
Y$ be the sets of abstract edges of $P$ and $R$, and $\ppi_X, \ppi_Y$
the projections to $\dq$.  Then we can cut and glue along subsets $I_X
\subset X, I_Y \subset Y$ equipped with a bijection $\phi: I_X
\rightarrow I_Y$ such that $\ppi_Y \circ \phi = * \circ \phi$, much in
the same way as the above: first, extend $\phi$ by $\phi^{-1}$ to
$I_Y$ to get an involution on $I_X \sqcup I_Y$. Then we take $X \sqcup
Y$, look at orbits of this under the map $f$ defined just as above,
and then define the map $\ppi'$ just as above (except that we need to
use $\ppi_Y$ instead of $\ppi_X$ on edges of $Y$), and
symmetric-multiply the result with any vertex idempotents appearing in
the original formulas for $P$ and $R$.

It is this latter operation which is what we precisely meant when we
spoke of ``cutting along edges and gluing the endpoints'' in the
definition of \eqref{act}.  In that case, we are summing over all
ordered choices of distinct elements $x_1, x_2, \ldots, x_m \in X$ and
$y_1, y_2, \ldots, y_m \in Y$, such that $\ppi_X(x_i) = e_i$ and
$\ppi_X(y_i) = e_i^*$.  Then we let $I_X = \{x_1, \ldots, x_m\}$ and
$I_Y = \{y_1, \ldots, y_m\}$ and $\phi(x_i) = y_i$, and perform
cutting and gluing (multiplying in some coefficient in $\k[\h]$).

Now that we have defined the action of \eqref{act}, we can extend
linearly over $\k$ to obtain the action of $\pi^n: \Sym L \o \Sym L
\rightarrow \Sym L$ for any $n$, and by linearity over $\k[\h]$, also
$e^{\frac{1}{2} \h \pi}: \Sym L[\h] \o \Sym L[\h] \rightarrow \Sym
L[\h]$. (Note that only polynomials in $\h$ are required since the
application of any differential operator of degree greater than the
total number of edges appearing in a given $P \o R$ is zero).

Now, we define $*_\h: \Sym L[\h] \o \Sym L[\h] \rightarrow \Sym L[\h]$ by
\begin{equation}
P *_\h R = e^{\frac{1}{2} \h \pi} (P \o R).
\end{equation}
This defines the necessary product which allows us to define $\Sym
L[\h]_\moy$.

We can describe this more directly as follows: again let $P, R$ be of
the form \eqref{pform} with sets of abstract edges $X, Y$, respectively,
and maps $\ppi_X: X \rightarrow \dq, \ppi_Y: Y \rightarrow \dq$. 
Then 
\begin{equation}
P *_\h R = \sum_{(I_X, I_Y, \phi)} \frac{\h^{\#(I_X)}}{2^{\#(I_X)}}
s(I_X, I_Y, \phi)
PR_{I_X, I_Y, \phi},
\end{equation}
where $(I_X, I_Y, \phi)$ is any triple of a subset $I_X \subset X, I_Y
\subset Y$ and a bijection $\phi: I_X \rightarrow I_Y$ satisfying
$\ppi_Y \circ \phi = * \circ \ppi_X$, and $PR_{I_X, I_Y, \phi}$ is the
result of cutting and gluing $P$ and $R$ along this triple as
described previously.  The sign $s(I_X, I_Y, \phi)$ is defined by
$s(I_X, I_Y, \phi) = (-1)^{\#(I_Y \cap \ppi_Y^{-1}(Q))}$. This follows
because $e^{\frac{1}{2} \h \pi} = \sum_{N \geq 0} \frac{\h^N}{2^N}
\frac{\pi^N}{N!}$, and each $\pi^N$ involves a sum over all cuttings and
gluings of $P$ and $R$ along $N$ edges counting each ordering and
multiplying in $-1$ for each time the $\frac{\partial}{\partial e}$
appears in the second component for $e \in Q$; dividing by $N!$ means
we don't count orderings of $I_X$ so that it is only over subsets that
we sum.

In general, elements $P, R \in \Sym L[\h]$ are linear combinations over
$\k[\h]$ of such collections of necklaces, so the element $P *_\h R$
is given by summing over each choice of necklace collections in $P$
and $R$, of the product of the coefficients of the two necklace
collections and the element described in the previous paragraph.  In
other words, we sum over all ways to take the product of terms from $P$
and $R$, not just by the usual product in $\Sym L[\h]$, but also by 
$\frac{\h^p}{2^p}$ times the ways in which we can cut out $p$ matching
edges from each term and join them together (while just multiplying
the $\k[\h]$-coefficients).

\subsection{Definition of the symmetrization map $\Phi_W$.} \label{phiws}
Now, we define $\Phi_W: \Sym L[\h] \rightarrow A_\h$.  To do this, we
need to define the notion of ``height assignments''.  Let's consider a
collection of necklaces $P$ of the form \eqref{pform}.  Let $X$ be the
set of abstract edges of $P$, say $\#(X) = N$.  Then, a \textsl{height
assignment} for $P$ is defined to be a bijection $H: X \rightarrow
\{1, 2, \ldots, N\}$.  
We have the
element $P_H \in A_\h$ obtained by assigning heights to the edges in
$X$ by $H$, that is
\begin{multline}
P_H = (a_{11}, H(1,1)) \cdots (a_{1 l_1},
H((1, l_1)) \& \cdots \\ \& (a_{k 1}, H(k, 1)) \cdots (a_{k l_k}, 
H(k, l_k)) \& v_1 \& v_2 \& \cdots \& v_q.
\end{multline}
Note that we could also think of $H$ as an element of $S_N$ with some
modifications to the formula.

The element $\Phi_W$ involves taking an average over all 
height assignments:
\begin{equation}
\Phi_W(P) = \frac{1}{N!} \sum_{H} P_H,
\end{equation}
where $H$ ranges over all height assignments.  Following is the alternative
description in terms of permutations $S_N$:
Let $\theta(i,j) = j + \sum_{p = 1}^{i-1} l_p$ so that $\theta(k,
 l_k) = N$. Then
\begin{multline}
\Phi_W(a_{11} \cdots a_{1 l_1} \& a_{21} \cdots
 a_{2 l_2} \& \cdots \& a_{k 1} \cdots a_{k l_k} \& v_1 \& v_2 \& \cdots \& v_q) \\ = \sum_{\sigma
 \in S_N} \frac{1}{N!} (a_{11}, \sigma(\theta(1,1))) \cdots (a_{1 l_1},
\sigma(\theta(1, l_1))) \& \cdots \\ \& (a_{k 1}, \sigma(\theta(k, 1))) \cdots
(a_{k l_k}, \sigma(\theta(k, l_k))) \& v_1 \& v_2 \& \cdots \& v_q.
\end{multline}

\subsection{Proof that $*_\h$ is obtained from $\Phi_W$.}\label{mpts}
Let's show that $\Phi_W$ makes the diagram \eqref{2d} commute. We know
that $\Phi_W$ is an isomorphism of free $\k[\h]$-modules (using PBW
for $A_\h$, or the fact that $\rho_\d$ is asymptotically injective and
the fact that the Weyl symmetrization map is an isomorphism on the
right-hand side of \eqref{2d}).  So, once we show commutativity of the
diagram, it will follow that $\Phi_W$ induces some multiplicative
structure on $\Sym L[\h]_\moy$ making the $\Phi_W$ an isomorphism of
$\k[\h]$-algebras.  We will then want to show that this structure is
the one we have just defined, i.e.~to show that $\Phi_W$ is a
homomorphism of rings using our $*_\h$ structure.

We need to show that $\rho_\d \circ \Phi_W = \phi_W \circ \tr$. This
follows immediately from the definitions, because $\rho_\d \circ
\Phi_W$ involves summing over the symmetrization of polynomials in
$(M_e)_{ij}, \frac{\partial}{\partial (M_e)_{ji}}, e \in Q$ where
$(M_e)_{ij}$ are the coordinate functions of the matrix corresponding
to the vertex $e$; also, $\tr$ takes an element of $\Sym L[\h]_\moy$
and gives the element of $\k[\h][Rep_\d(\dq)]$ corresponding to the
trace of the (cyclic noncommutative) polynomial, which upon
substituting $(M_{e^*})_{ij} \mapsto -h \frac{\partial}{\partial (M_e)_{ji}}$
and symmetrizing (which we needed to do for this to be well-defined,
since the $(M_{e^*})_{ij}, (M_e)_{ij}$ commuted), gives the same element.

Next, let us show that the ring structure obtained from $\Phi_W$, making
$\Phi_W$ an isomorphism of rings, is exactly the product $*_\h$ we
have described in detail.
\begin{equation}
\label{pqe}
\Phi_W(P *_\h R) = \Phi_W(P) \Phi_W(R).
\end{equation}

Now we prove \eqref{pqe}.  Let's take $P = P_1 \& P_2 \& \cdots \&
P_n$, as before, to be a collection of necklaces, and similarly for $R
= R_1 \& R_2 \& \cdots \& R_m$.  (We can forget about the idempotents
such as in \eqref{pform}, since they won't affect what we have to
prove.)  Let $X$ be the set of abstract edges of $P$ and $Y$ the set
of abstract edges of $R$.  We will use $H_P: X \rightarrow \{1,
\ldots, |X|\}$ to denote a height assignment for $P$ and $H_R: Y
\rightarrow \{1, \ldots, |Y|\}$ to denote a height assignment for $R$.
Let us say that a height assignment $H_{PR}: X \sqcup Y \rightarrow
\{1, \ldots, |X|+|Y|\}$ \textsl{extends} height assignments $H_P, H_R$
if $H_{PR}$ restricted to $P$ is equivalent to $H_P$ and $H_{PR}$
restricted to $R$ is equivalent to $H_R$. In other words, $H_{PR}(x_1)
< H_{PR}(x_2)$ iff $H_P(x_1) < H_P(x_2)$ for all $x_1, x_2 \in X$, and
similarly $H_{PR}(y_1) < H_{PR}(y_2)$ iff $H_R(y_1) < H_R(y_2)$ for
all $y_1, y_2 \in Y$.

Now, we know that
\begin{equation}
\Phi_W(P *_\h R) - \Phi_W(PR) = \sum_{N = 1}^\infty \frac{\h^N}{2^N} 
\Phi_W(\frac{\pi^N}{N!} (P \o R)),
\end{equation} 
and also that
\begin{multline}
\Phi_W(P) \Phi_W(R) - \Phi_W(PR) \\ = \frac{1}{(|X|+|Y|)!} 
\sum_{H_P, H_R} \sum_{H_{PR} \text{\ extending\ }H_P, H_R} (P_{H_P} R_{H_R}
- PR_{H_{PR}}).
\end{multline}

We are left to show, using the relations which define $A_\h$, that
\begin{equation} \label{ltsiso}
\sum_{N = 1}^\infty \frac{\h^N}{2^N} \Phi_W(\frac{\pi^N}{N!} (P \o R)) = 
\sum_{H_{PR} \text{\ extending\ }H_P, H_R} (P_{H_P} R_{H_R}
- PR_{H_{PR}})
\end{equation}

To prove this, let us fix a particular $H_P, H_R$, and $H_{PR}$, and 
expand $P_{H_P} R_{H_R} - PR_{H_{PR}}$ using
the relations that define $A_\h$.  We do this by expressing this as a sum
of commutators obtained by commuting a single edge coming from $R$
with a single edge coming from $P$.  We get
\begin{equation}
P_{H_P} R_{H_R} - PR_{H_{PR}} = \sum_{\underset{\ppi_X(x) = \ppi_Y(y)^*}{x \in X, y \in Y \text{\ such that\ 
} H_P(x) > H_R(y),} } [\ppi_X(x), \ppi_Y(y)] \h PR'_{x,y},
\end{equation}
where $PR'_{x,y}$ corresponds to taking $PR$, deleting $x$ and $y$ and
joining the endpoints, and using the height assignment which is
(equivalent to the choice of heights) identical to $H_P$ on elements
$x' \in X$ such that $H_P(x') < H_P(x)$, and equal to $H_P(x) +
H_{PR}(z)$ for all other $z \in X \sqcap Y \setminus \{x,y\}$.  Here
we say ``equivalent to the choice of heights'' in parentheses to mean
that the given assignment won't be an assignment to $\{1, \ldots,
|X|+|Y|-2\}$ as we strictly defined height assignments, but we could
extend the definition of height assignments to include any injective
map to $\Z$, and say that two are equivalent if the ordering is the
same ($H \equiv H'$ if $H(z) < H(z')$ iff $H'(z) < H'(z')$); so really
we are looking for the height assignment mapping to $\{1, \ldots,
|X|+|Y|-2\}$ which is equivalent to the assignment we described. Also
note here that $[e,e^*] = 1$ if $e \in Q$ and $-1$ if $e^* \in Q$.

By applying the relations repeatedly we get that
\begin{multline} \label{hpre}
P_{H_P} R_{H_R} - PR_{H_{PR}} \\ = \sum_{\underset{\text{\ such that\ }
 H_P(x_i) > H_R(y_i), \ppi_X(x_i)
= \ppi_Y(y_i)^*}{x_1, \ldots, x_k \in X, y_1,
\ldots, y_k \in Y } } [\ppi_X(x_1), \ppi_Y(y_1)] \\ \cdots [\ppi_X(x_k),
\ppi_Y(y_k)] \h^k PR''_{(x_1, y_1), \ldots, (x_k, y_k)},
\end{multline}
where $PR''_{(x_1, y_1), \ldots, (x_k, y_k)}$ involves taking $PR$ and
deleting the pairs of edges and gluing at their respective endpoints;
and this time assigning heights by restricting $H_{PR}$ to $X \sqcup Y
\setminus \{x_1, \ldots, x_k, y_1, \ldots, y_k\}$ (and changing to an
equivalent height assignment which has image $\{1, \ldots,
|X|+|Y|-2k\}$).

Now, let's look at the sum again (no longer fixing $H_P, H_R$, and
$H_{PR}$).  We see that for any given choice of pairs $(x_1, y_1),
\ldots, (x_k, y_k)$ with $\ppi_X(x_i) = \ppi_Y(y_i)^*$, the summands
that involve deleting these pairs and gluing their endpoints are the
same in number for each choice of height assignment for the deleted
pairs.  The coefficient for each height is just
$\frac{\h^k}{(|X|+|Y|)!}$ times the number of height assignments
$H_{PR}$ that restrict to the given height assignment, and also
satisfy $H_{PR}(x_i) > H_{PR}(y_i)$ for all $1 \leq i \leq k$.  In
other words, this is $\h^k$ times the probability of picking a height
assignment randomly of $PR$ that has $x_i$ greater in height than
$y_i$ for all $i$, and is identical with the given height assignment
on all $x, y \notin \{x_1, \ldots, x_n, y_1, \ldots, y_n\}$.  So we
get that the coefficient is $\frac{\h^k}{2^k (|X|+|Y|-2k)!}$. 

But this is exactly what we would expect, desiring that \eqref{ltsiso}
hold. That is because the left-hand side, as described previously in
our discussion of $\frac{\pi^N}{2^N}$, just involves summing over all
$N$ of $\frac{\h^k}{2^k}$ times $\Phi_W$ of the collection of
necklaces described for each choice of pairs $(x_1, y_1), \ldots,
(x_N, y_N)$ with some choice of sign; and then $\Phi_W$ just sums over
$\frac{1}{(|X|+|Y|-2N)!}$ times each choice of height assignment for
this collection of necklaces. The sign choice just matches exactly with
the sign $\prod_{i} [\ppi_X(x_i), \ppi_Y(y_i)]$ appearing in \eqref{hpre},
since each commutator is $-1$ just in the case that $\ppi_Y(y_i) \in Q$.

This proves \eqref{ltsiso} and hence that $\Phi_W$ is an isomorphism
of $\k[\h]$-algebras, using $*_\h$ as the ring structure on $\Sym L[\h]_\moy$.

\subsection{Associativity of $*_\h$.} \label{ass}
Although we already know from commutativity of the diagram and associativity
of $A_\h$ that $*_\h$ is associative, it is easy to prove directly. We
prove
\begin{equation} \label{tpass}
(P *_\h R) *_\h S = P *_\h (R *_\h S)
\end{equation}
where $P, R$, and $S$ are collections of necklaces of the form \eqref{pform}
(with different indices).

First we describe the left-hand side of \eqref{tpass}
Let $X, Y$, and $Z$ be the sets of abstract edges
of $P, R,$ and $S$, and let $\ppi_X: X \rightarrow \dq, \ppi_Y: Y
\rightarrow \dq$, and $\ppi_Z: Z \rightarrow \dq$ be the projections
from occurrences of edges to edges of $\dq$.

We sum over all sets of pairs $\{(x_1, y_1), \ldots, (x_N, y_N)\}
 \subset X \times Y$, such that $y_i = x_i^*$ for each $i$ (and we
 assume that the $x_i$ and the $y_i$ are all distinct).  Summing over
 $\frac{\h^N}{2^N}$ times the necklaces we get by cutting out these
 pairs of edges and gluing their endpoints, we get $P *_\h R$ as
 described in the previous section.

To get $(P *_\h R) *_\h S$, we will first be summing over choices of
pairs $\{(x_1, y_1), \ldots (x_N, y_N)\}$, and then over pairs
$\{(w_1, z_1), \ldots, (w_M, z_M)\} \subset W \times Z$, where $W = (X
\setminus \{x_1, x_2, \ldots, x_N\}) \sqcup (Y \setminus \{y_1, y_2,
\ldots, y_N\})$, and performing a similar operation.  We can also
describe this as summing over pairs $(x_1, y_1), \ldots, (x_N, y_N),
(x_{N+1}, z_1), (x_{N+2}, z_2), \ldots, (x_{N+M_1}, z_{M_1}),$ \\
$(y_{N+1}, z_{M_1+1}), (y_{N+2}, z_{M_1+2}), \ldots, (y_{N+M_2},
z_{M_1+M_2})$, again where all $x_i, y_i,$ and $z_i$ are distinct, and
in each pair, one edge is the reverse of the other.  This description,
along with signs and coefficients, is exactly the same as what we
could obtain in the same way from $P *_\h (R *_\h S)$, proving
associativity.

\section{The Moyal coproduct} 
\subsection{Definition of $\Delta_\h$.}\label{cps}
There is a nice formula for the coproduct on $\Sym L[\h]_\moy$ compatible
with the the $*_\h$ product.  The formula is actually surprisingly
similar to the one for $*_\h$.  We will be giving the coproduct which
makes the diagram \eqref{2d} consist of coalgebra homomorphisms
(namely, the maps $\tr$ and $\Phi_W$ involving $\Sym L[\h]_\moy$); the map
$\Phi_W$ will then be an isomorphism of bialgebras.  The coproduct can
be described as follows: We need to define an operator $e^{\frac{1}{2}
\h \pi}: \Sym L[\h]_\moy \rightarrow \Sym L[\h]_\moy \o \Sym L[\h]_\moy$.  To do
this, we set
\begin{equation}
\pi = \sum_{e \in Q} \frac{\partial}{\partial e} \frac{\partial}{\partial e^*}
\end{equation}
and we define operators
\begin{multline} \label{bco}
\frac{\partial}{\partial e_1} \frac{\partial}{\partial e_1^*} \frac{\partial}{\partial e_2} \frac{\partial}{\partial e_2^*} \cdots \frac{\partial}{\partial e_N}\frac{\partial}{\partial e_N^*}: \\ \Sym L[\h]_\moy \rightarrow \Sym L[\h]_\moy \o \Sym L[\h]_\moy.
\end{multline}
The operator \eqref{bco} acts as follows: Taking a collection of
necklaces $P = P_1 \& P_2 \& \cdots \& P_n$ \\ $\& v_1 \& v_2 \& \cdots \& v_q
\in \Sym L[\h]_\moy$, where
each $P_i \in L$ is a cyclic monomial (i.e.~a necklace), let $X$ be
the set of abstract edges of $P$ and $\ppi_X: X \rightarrow \dq$ the
projection (cf.  Section \ref{mps}).  Then we sum over all choices of
pairs $(x_1, y_1), \ldots, (x_N, y_N)$ such that the $x_i$ and $y_i$
are all distinct (the set $\{x_1, y_1, \ldots, x_N, y_N\}$ has $2N$
elements), and $\ppi_X(x_i) = \ppi_X(y_i)^*$ for all $i$. We delete the
edges and glue the endpoints, obtaining another collection of
necklaces.  More precisely, the cutting and gluing is done as
described in the previous section, for $I = \{x_1, y_1, x_2, y_2,
\ldots, x_N, y_N\}$ and $\phi(x_i) = y_i$ for all $i$.  Now, the only
difficult part is figuring out what components to assign to each
necklace (the first or second), and what sign to attach to each
choice.

We sum over all component assignments of the resulting chain of
necklaces: suppose that the above procedure yields the collection $R_1
\& R_2 \& \cdots \& R_m \in \Sym L[\h]_\moy$ (this includes the original
idempotents $v_1, v_2, \ldots, v_q$); then the contribution to the
result of \eqref{bco} applied to $P$ is the following:

\begin{equation} \label{cass}
\sum_{\mathbf{c} \in \{1,2\}^m} s(\mathbf c, I, \phi)
 R_1^{c_1} \& R_2^{c_2} \& \cdots \& R_m^{c_m},
\end{equation}
where $R_i^{c_i}$ denotes $R_i \o 1$ if $c_i = 1$ and $1 \o R_i$ if
$c_i = 2$, and the symmetric product in $\Sym L[\h]_\moy \o \Sym L[\h]_\moy$
is the expected $(X \o Y) \& (X' \o Y') = (X \& X') \o (Y \& Y')$,
with $1 \& X = X \& 1 = X$ for all $X$.  The term $s(\mathbf c, I, \phi)$ 
is a sign which is determined as follows: 
\begin{equation}
s(\mathbf c, I, \phi) = s_1 s_2 \cdots s_n,
\end{equation}
where $s_i = 1$ if the component assigned to the start of
arrow $x_i$ is $1$ and the component assigned to the target of arrow
$x_i$ is $2$; $s_i = -1$ if the component assigned to the start of
arrow $x_i$ is $2$ and the component assigned to the target of arrow
$x_i$ is $1$; and $s_i = 0$ if the start and target are assigned the
same component.  

Let's more precisely define what it means to say ``the component
assigned to the start/target of an arrow'' which is deleted from $P$.
What we mean by this is simply the orbit of the arrow $x_i$ in $X$
under $f$.  Each orbit corresponds to one of the $R_i$.  So, there is
a map $g: X \rightarrow \{1, 2, \ldots, m\}$, which corresponds to
which $R_i$ the ``start'' of each edge is assigned to.  The
``targets'' are the same as the ``starts'' of the next edge, so that
$g(x+1)$ gives the component that the ``target'' of $x$ is assigned
to. Here the ``$+1$'' operation is once again $(i,j)+1=(i,j+1)$ mod
$l_i$, or in other words, $x+1$ is the edge succeeding $x$.

We then have that
\begin{equation} \label{scc}
s_i = \begin{cases} 1 & c_{g(x_i)} < c_{g(x_i)+1}, \\
                    0 & c_{g(x_i)} = c_{g(x_i)+1}, \\
                   -1 & c_{g(x_i)} > c_{g(x_i)+1}.
\end{cases}
\end{equation}

This assignment of signs has a combinatorial flavor because it is
essentially what the ``colorings'' of \cite{S} reduce to.  There does
not seem to be a way to avoid this complication in choosing signs,
because the sign is what prevents the coproduct from being
cocommutative.

As before, we extend linearly to powers $\pi^N$ and to $e^{\frac{1}{2} \h \pi}$.
Then, the coproduct is given by
\begin{equation}
\Delta_\h = e^{\frac{1}{2} \h \pi}: \Sym L[\h]_\moy \rightarrow \Sym L[\h]_\moy \o \
\Sym L[\h]_\moy,
\end{equation}
and as before we can describe this action on our element $P$ as
\begin{equation}
\Delta_\h(P) = \sum_{(I, \phi)} \frac{\h^{\#(I)/2}}{2^{\#(I)/2}} P_{I, \phi}, \\ 
\end{equation}
where the sum is over all $I \subset X$ with involution $\phi$ such
that $\ppi_X \circ \phi = * \circ \ppi_X$, and the element $P_{I, \phi}$
is given from the result of the cuttings and gluings by summing over
component assignments as described in \eqref{cass}.

\subsection{Proof that $\Delta_\h$ is obtained from $\Phi_W$.} \label{cpts}
Let's prove that this coproduct $\Delta_\h$ makes the diagram \eqref{2d}
consist of coalgebra homomorphisms.  It suffices to prove that
$\Phi_W$ is a coalgebra homomorphism.

Take an element $P$ of the form \eqref{pform} with set of abstract
edges $X$ and projection $\ppi_X: X \rightarrow \dq$.  Now, let's
consider what the element $\Delta(\Phi_W(P))$ is in $A$. We know that
for each height assignment $H_P$ of $P$, $\Delta(P_{H_P})$ involves
summing over all pairs $(I, \phi)$ with $I \subset X$ and $\phi: I
\rightarrow I$ an involution satisfying $\ppi_X \circ \phi = * \circ
\ppi_X$, cutting and gluing as before. Then we sum over all component
assignments such that if $x, y \in I$ with $\phi(x) = y$, and the
heights satisfy $H(x) < H(y)$, then the component assigned to the
start of $x$ is $1$ and the component assigned to the target of $x$ is
$2$. When the components cannot be assigned in this way, this pair
$(I, \phi)$ cannot be used.  These notions are all explained more
precisely in the preceding section.

Then we multiply in a sign $s(I, \phi, H)$ and a power of $\h$ determined
as follows: for each pair $x, y \in I$ with $\phi(x) = y, H(x) <
H(y)$, we multiply a $+1$ if $x \in Q, y \in Q^*$ and a $-1$ if $x \in
Q^*, y \in Q$.  We also multiply in $\h^{\#(I)/2}$ (note: this power
of $\h$ is different from the one in \cite{S} simply because we are
describing the structure for $A_\h$, not $A$: it is easy to see in
general how the relations for the algebra and the formula for
coproduct change if we introduce a new formal parameter $\hbar$ into 
(3.3) of \cite{S}).

So we find that $\Delta(P_H)$ is just a sum over cuttings and gluings,
and over component choices $\mathbf c$ compatible with the heights;
our sign choice satisfies $s(I, \phi, H) = s(\mathbf c, I, \phi)$,
where $I = \{x_1, y_2, \ldots, x_m, y_m\}$, and for all $i$, $x_i \in
Q$ and $\phi(x_i) = y_i$; finally, we multiply in $\h^m$ for cuttings
and gluings involving $\#(I)=2m$.

Hence, $\Delta(\Phi_W(P))$ is just given by a sum over all cuttings and
gluings $(I, \phi)$ together with component choice $\mathbf c$, 
multiplying in $\h^{\#(I)/2}$, the sign $s(\mathbf c, I, \phi)$, and
the coefficient $\frac{1}{\#(P)!}$ where $\#(P)$ is the number of edges
in $P$, i.e.~the total number of height assignments. 

Each summand in $\Delta(\Phi_W(P))$ is clearly given by a height
assignment of the term in $\Delta_\h(P)$ corresponding to the same $(I,
\phi, \mathbf c)$.  For each term in $\Delta_\h(P)$, the coefficients of
the height-assigned terms in $\Delta(\Phi_W(P))$ are all the same.  So
we see that $\Delta(\Phi_W(P)) = (\Phi_W \o \Phi_W)(P')$, for some $P'
\in \Sym L[\h]_\moy \o \Sym L[\h]_\moy$, where $\Phi_W \o \Phi_W (P \o R) =
\Phi_W(P) \o \Phi_W(R)$.

The element $P'$ can be computed just as we were computing $\Delta(\Phi_W(P))$, but instead of multiplying in $\frac{1}{\#(P)!}$, we need to multiply by the 
fraction of all height choices compatible with this
component choice.  But clearly, each pair $x, y \in I, \phi(x) = y$
induces a single restriction on the choice of heights, namely that
$H(x) < H(y)$ if the component assigned to the start of $x$ is $1$
and the component assigned to the target of $x$ is $2$, and $H(y) > H(x)$
if the opposite is true (the start of $x$ is assigned component $2$ and
the target assigned $1$).  Note that the component assigned to the start
and target of $x$ cannot be the same for there to exist any compatible
height choices.

We see then that, provided a compatible height choice exists, we have
$\#(I)/2$ restrictions, each of which occur with independent probabilities
$\frac{1}{2}$. Hence the coefficient is just $\frac{1}{2^{\#(I)/2}}$. This
shows that $P' = \Delta_\h(P)$, proving that $\Phi_W$ is a coalgebra homomorphism
and hence an isomorphism of bialgebras. (In fact we have now proved that
$(\Sym L[\h]_\moy, *_\h, \Delta_\h)$ is in fact a bialgebra).

\subsection{Coassociativity of $\Delta_\h$.} \label{casss}
Using the coassociativity of $A_\h$ from \cite{S}, we already know from
the fact that $\Phi_W$ is an isomorphism that the product $\Delta_\h$ is
coassociative, but it is not difficult to prove directly. We do that here
by proving
\begin{equation} \label{coasse}
(1 \o \Delta_\h) \Delta_\h (P) = (\Delta_\h \o 1) \Delta_\h (P),
\end{equation}
where once again $P$ is of the form \eqref{pform}.

The left-hand side can be described by summing over choices of cutting
pairs and components $(I, \phi, \mathbf c)$ for $P$, and then cutting
pairs and components for the first component of the result, $(I',
\phi', \mathbf {c'})$, and gluing, assigning the components, etc., and
multiplying by a sign and power of $\frac{\h}{2}$.  We see that this
is the same as choosing just once the triple $(I'', \phi'', \mathbf
{c''})$, where $\mathbf {c''}$ assigns each necklace to one of three
components, $1, 2,$ or $3$, $I'' = I \cup I'$, and $\phi'' |_I = \phi,
\phi''|_{I'} = \phi'$. Then we can cut and glue just one time to get a
tensor in $\Sym L[\h]_\moy^{\o 3}$; the sign and power of $\frac{\h}{2}$
can be determined by using \eqref{scc} where now the two sides of the
inequality have values in $\{1,2,3\}$.

For the same reason, the right-hand side of \eqref{coasse} can be
described in the preceding way, proving \eqref{coasse} and hence
coassociativity.

\subsection{The bialgebra condition for $\Sym L[\h]_\moy$.} \label{bas}
Using the fact that $A_\h$ is a bialgebra (proved in \cite{S}), we
know that $\Sym L[\h]_\moy$ is a bialgebra. But it is not difficult to
prove directly, which we do in this section.

We need to show, for collections of necklaces 
$P = P_1 P_2 \cdots P_n, R = R_1 R_2 \cdots R_m$ of the form \eqref{pform},
that
\begin{equation} \label{bac}
\Delta(P *_\h R) = \Delta(P) *_\h \Delta(R),
\end{equation}
where we use the notation $(A \o B) *_\h (C \o D) = (A *_\h B) \o (C *_\h D)$.

First, define $X$ to be the set of abstract edges of $P$ and $Y$ to be
the set of abstract edges for $R$. Define the projections $\ppi_X: X
\rightarrow \dq, \ppi_Y: Y \rightarrow \dq$. Now, let us first take a
closer look at the right-hand side of \eqref{bac}.  We can expand it
by the following sum over pairs. We first pick $I_X \subset X, I_Y
\subset Y,$ and involutions $\phi_X: I_X \rightarrow I_X, \phi_Y: I_Y
\rightarrow I_Y$ such that $\ppi_X \circ \phi_X = * \circ \ppi_X, \ppi_Y
\circ \phi_Y = * \circ \ppi_Y$.  Pick component choices $\mathbf{c}$
for the result of cutting and gluing $P$ along $(I_X, \phi_X)$, and
$\mathbf{c'}$ for the result of cutting and gluing $R$ along $(I_Y,
\phi_Y)$.  As before, we define signs $s(I_X, \phi_X, \mathbf c),
s(I_Y, \phi_Y, \mathbf {c'})$. For example, $s(I_X, \phi_X, \mathbf c)$
is defined by multiplying in all the $\pm 1$ or $0$ contributions from
each $x \in I_X$ such that $\ppi_X(x) \in Q$, according to \eqref{scc}.
Next, we cut and glue both $P$ and $R$ by $(I_X, \phi_X, \mathbf c)$
and $(I_Y, \phi_Y, \mathbf {c'})$, respectively, and multiply the first
by $s(I_X, \phi_X, \mathbf c)\frac{\h^{\#(I_X)/2}}{2^{\#(I_X)/2}}$ and
the second by $s(I_Y, \phi_Y, \mathbf {c'}) \frac{\h^{\#(I_Y)/2}}
{2^{\#(I_Y)/2}}$ to obtain elements $P', R' \in \Sym L[\h]_\moy$.
This will include all the summands we need for the coproducts of $P$ and
$R$, respectively.

For each such summand, we need to take care of contributions from
multiplying these together. So, we need to pick $J_X \subset X
\setminus I_X, J_Y \subset Y \setminus I_Y$, and a bijection $\psi:
J_X \rightarrow J_Y$ such that $\ppi_Y \circ \psi = * \circ \ppi_X$ and
also the extra condition that $\psi$ preserves components: that is, if
$\psi(x) = y$ and $x, y$ live in necklaces assigned components $c_i,
c_j'$, respectively, then $c_i = c_j'$.

To be more precise, the cutting and gluing $P \mapsto T_1 \& T_2 \&
\cdots \& T_p$ along $(I_X, \phi_X)$ induces a map $\mu_X: X \setminus
I_X \rightarrow \{1, 2, \ldots, p\}$ depending on which necklace each
edge not cut out ends up in.  So each edge $x \in X \setminus I_X$ is
assigned a component $c_{\mu_X(x)}$.  Similarly we can define $\mu_Y$.
The condition above is that $\psi(x) = y$ implies that $c_{\mu_X(x)} =
c'_{\mu_Y(y)}$.

Given each such choice of $(I_X, \phi_X, \mathbf c), (I_Y, \phi_Y,
\mathbf {c'}),$ and $(J_X, J_Y, \psi)$, we get the following
contribution to the expression $\Delta(P) *_\h \Delta(R)$: We cut and
glue $P'$ and $R'$ together along $(J_X, J_Y, \psi)$, and multiply
in the sign $s(J_X, J_Y, \psi) = (-1)^{\#(J_Y \cap \ppi_Y^{-1}(Q))}$
and the coefficient $\frac{\h^{\#(J_X)/2}}{2^{\#(J_X)/2}}$.  We sum
this contribution over all such triples of triples to get the right-hand
side of \eqref{bac}.  

Now, let's compare this with the left-hand side of \eqref{bac}.  Let the
map $\ppi: X \sqcup Y \rightarrow \dq$ be given by $\ppi(x) = \ppi_X(x),
\ppi(y) = \ppi_Y(y)$ for $x \in X, y \in Y$. Then, the left-hand side
involves first a sum over $(J_X, J_Y, \psi)$ with $J_X \subset X, J_Y
\subset Y,$ and $\psi: J_X \rightarrow J_Y$ a bijection such that $\ppi_Y
\circ \psi = * \circ \ppi_X$.  Then we sum over $(I, \phi, \mathbf {c''})$
such that $I \subset (X \setminus J_X) \sqcup (Y \setminus J_Y)$ and
$\phi: I \rightarrow I$ is an involution satisfying $\ppi \circ \phi = *
\circ \ppi$, and $\mathbf {c''}$ is a component choice of $PR_{J_X, J_Y,
\psi}$, the result of
cutting and gluing $P$ and $R$ along $(J_X, J_Y, \psi)$ and then along
$(I, \phi)$.  For each such choice of triples, we multiply a coefficient
of $\frac{\h^{\#(J_X) + \#(I)/2}}{2^{\#(J_X) + \#(I)/2}}$, and a sign of
$(-1)^{\#(J_Y \cap \ppi_Y^{-1}(Q))} s(I, \phi, \mathbf {c''})$ where the
$s(I, \phi, \mathbf {c''})$ is calculated just as in the definition of 
$\Delta(PR_{J_X, J_Y, \psi})$.

If $\phi(I \cap X) \subset X$ and $\phi(I \cap Y) \subset Y$, then
the summand thus obtained will be identical with the summand
corresponding to $(I \cap X, \phi|_X, \mathbf c), (I \cap Y, \phi|_Y,
\mathbf {c'}), (J_X, J_Y, \psi)$ for the choice of $\mathbf c, \mathbf
{c'}$ such that $c''_{\mu_{X \sqcup Y \setminus (J_X \cup J_Y)}(z)}$
equals $c_{\mu_X(z)}$ if $z \in X$ and $c'_{\mu_Y(z)}$ if $z \in Y$,
where $\mu_{X \sqcup Y \setminus (J_X \cup J_Y)}$ is defined just as
$\mu_X, \mu_Y$ were, in the context of cutting and gluing $PR_{J_X, J_Y,
\psi}$ along $(I, \phi)$.  The power of $\frac{\h}{2}$ will clearly be
the same.  The sign will also be the same, since $s(I, \phi, \mathbf
{c''}) = s(I \cap X, \phi|_X, \mathbf c) s(I \cap Y, \phi|_Y, \mathbf
{c'})$ in this case.

All that remains is to show that all summands from the left-hand side
not of this form cancel. Summands which are not of the above form must
either include $x \in I \cap X, y \in I \cap Y$ such that $\phi(x) = y$,
or else must have some component choice such that $c''_i \neq c''_j$
even though the necklaces $i$ and $j$ would be joined if we had omitted
some $x$ from $J_X$ and $\psi(x)$ from $J_Y$.  The latter comes from the
fact that $c''_i = c''_j$ is exactly what is required for $\mathbf{c''}$
to be compatible with some $\mathbf{c}, \mathbf{c'}$ such that
$c_{\mu_X(x)} = c'_{\mu_Y(y)}$.

Let us make the definitions
\begin{multline}
J_X' =  
\{x \in J_X \mid c''_i \neq c''_j, \\ \text{\ where $i$ and $j$ would
be joined by omitting $x, \phi(x)$ from $I$}\},
\end{multline}
\begin{equation}
I_X' = \{x \in I \cap X \mid \phi(x) \in Y\}.
\end{equation}
For each $x \in J_X'$,
we can obtain a similar summand by removing $x$ from $J_X$ and $y = \psi(x)$
from $J_Y$, and adding $x, y$ to $I$, setting $\phi(x) = y, \phi(y) =
x$: so $x$ ends up in $I_X'$. 
We get the same resulting necklaces and can consider the
$\mathbf{c''}$ which makes the same assignments to the corresponding
necklaces (where necklaces correspond if they come from the same edges
in $X$ or vertex idempotents in the expression for $P$).  The only
change is perhaps a change of sign; the sign changes iff $c''_{g(x)} >
c''_{g(x)+1}$ where $g$ is defined as in \eqref{scc} (and $g(x) +1$ is
the edge following $x$ in $P$): we are saying that
the sign changes iff the component assigned to the start of $x$ is $2$
and the component assigned to the target of $x$ is $1$.

Similarly, for each such summand, we can perform the operation of
removing an $x \in I_X'$ and $y = \phi(x) \in (I \cap Y)$, and
adding $x$ to $J_X$ and $y$ to $J_Y$ in such a way that the component
assignments remain the same: in this case, $x$ ends up in $J_X'$. Again,
we get a sign change just in the event that $c''_{g(x)} > c''_{g(x)+1}$.

So, if we sum up all summands which can be obtained from each other by
applying the above two operations, we will get zero unless all sign
changes in the above two paragraphs are positive.  But, this cannot
happen if $J_X'$ or $I \cap X$ was originally nonempty for the following
reason:

\begin{equation} \label{fbap}
0 = \sum_{x \in X} (c''_{g(x)} - c''_{g(x)+1}) = \sum_{x \in I_X'
\cup J'_X} c''_{g(x)} - c''_{g(x)+1},
\end{equation}

so that one of the summands on the right-hand side must be negative
(since $J'_X \cup I_X' \neq \emptyset$ shows that the last summand is
nonempty, and each summand is $\pm 1$ in that last summand).  The
justification for passing from the first to the second summation in
\eqref{fbap} is that the only nonzero terms we have eliminated in doing
so are those that correspond to $x \in J_X \cup I \cap X$ such that $x$
is paired by $\phi$ or $\psi$ with another $x' \in X$, so that the first
summation includes both $c''_{g(x)} - c''_{g(x) + 1}$ and $c''_{g(x')} -
c''_{g(x') + 1}$, which cancel pairwise.

We have proven that all summands in the left-hand side of \eqref{bac}
either correspond in a bijective fashion with a summand from the
right-hand side, or else lie in a set of summands with nonempty $I_X'$
and $J_X'$, whose contributions cancel.  The proof of \eqref{bac} is
finished, so that $\Sym L[\h]_\moy$ is a bialgebra.

\section{The antipode} \label{as}
Using $\Phi_W$ and the formula for the antipode in \cite{S}, it
immediately follows that our antipode $S: \Sym L[\h]_\moy$
is given by the formula
\begin{equation} \label{man}
S(P_1 \& P_2 \& \cdots \& P_m) = (-1)^m P_1 \& P_2 \& \cdots \& P_m,
\end{equation}
where each $P_i \in L$ is a necklace (i.e.~a cyclic monomial or vertex
idempotent).  It is immediate that $S^2 = \Id$.  Indeed, $S$ is
diagonalizable with eigenvalues $\pm 1$ and eigenvectors which are
collections of necklaces of the form \eqref{pform}.

Unfortunately, a direct proof that \eqref{man} is the antipode for $\Sym
L[\h]_\moy$ turned out to be too difficult.  The authors are
interested in any good proof of this fact from purely the Moyal point of
view.

\bibliography{moyal}
\bibliographystyle{amsalpha}

\footnotesize{
{\bf V.G.}: Department of Mathematics, University of Chicago, 
Chicago IL
60637, USA;\\ 
\hphantom{x}\quad\, {\tt ginzburg@math.uchicago.edu}}
\smallskip 

\footnotesize{
{\bf T.S.}: Department of Mathematics, University of Chicago, 
Chicago IL
60637, USA;\\ 
\hphantom{x}\quad\, {\tt trasched@math.uchicago.edu}}

\end{document}